\documentclass[a4paper,11pt]{article} 
\usepackage[utf8]{inputenc}
\usepackage[T1]{fontenc}
\usepackage{amsmath,amsthm,amsfonts,amssymb,euscript,stmaryrd,graphicx,enumitem,yhmath,mathrsfs}
\usepackage[toc,page]{appendix}
\usepackage{algorithm,algorithmic}
\usepackage{fullpage}
\usepackage{color,psfrag}
\usepackage{bbm}

\usepackage{dsfont}
\usepackage[babel=true]{csquotes} 
\usepackage{lmodern}
\usepackage{bibentry}

\usepackage[french,english]{babel}
\usepackage{ifthen}
\usepackage{array}
\usepackage{hyperref} 
\usepackage{mathrsfs}

\usepackage{yfonts}

\def\ben{\begin{enumerate}}
\def\een{\end{enumerate}}

\newcommand{\bit} {\begin{itemize} }
\newcommand{\eit} {\end{itemize} }

\newcommand{\ds}{\displaystyle}

\newcommand{\fdd}{\stackrel{f.d.d.}{\longrightarrow}}
\newcommand{\pn}{\stackrel{\mathbb{P}}{\longrightarrow}}
\newcommand{\weak}{\stackrel{d}{\longrightarrow}}

\newcommand{\lin}{\underset{{n \to +\infty}}{\underline{\lim}}}
\newcommand{\lsp}{\ds \underset{{n \to +\infty}}{\overline{\lim}}}

\newcommand{\sifbm}{\mathbf{B}}

\def\i1{\mathbf{1}}  

\newcommand{\N}{\mathbb{N}}

\newcommand{\R}{\mathbb{R}}

\newcommand{\bN}{\mathbb{N}}

{\end{array}}

\newcommand{\footnoteremember}[2]{
 \footnote{#2}
 \newcounter{#1}
 \setcounter{#1}{\value{footnote}}
}

\author
{        {Joachim \textsc{Lebovits}\footnote{Laboratoire Analyse, Géométrie et Applications, C.N.R.S. (UMR 7539), Université Paris 13, Sorbonne Paris Cité, 99 avenue Jean-Baptiste Clément
93430, Villetaneuse, France. Email address: \url{jolebovits@gmail.com}.}}
    \and
        Mark \textsc{Podolskij}\footnoteremember{footnoteINRIA}{Department of Mathematics, University of Aarhus, Ny Munkegade 118, 8000 Aarhus C, Denmark. Email address: \url{mpodolskij@math.au.dk}}}

\title{Estimation of the global regularity of a multifractional Brownian motion}

\newtheorem{theo}{Theorem}[section]

\newtheorem{prop}[theo]{Proposition}

\newtheorem{rem}{Remark}

\newtheorem{defi}{Definition}

\newcommand{\E}{\mathbb{E}}

\newcommand{\cM}{\mathcal{M}}

\newcommand{\uh}{\mathscr{U}^{\hspace{-0.2ex}h}_{n}}

\newcommand{\card}{\operatorname{card}}

\def\h1{\hspace{0.1cm}}

\def\keywordname{{\bf Keywords:} consistency, Hurst parameter, multifractional Brownian motion, power variation}

\newcommand{\keywords}[1]{\par\addvspace\baselineskip\noindent\keywordname\enspace\ignorespaces#1}

\makeatletter
\renewcommand\theequation{\thesection.\arabic{equation}}
\@addtoreset{equation}{section}
\makeatother

\graphicspath{{./}{figures/}}

\begin{document}

\maketitle

			\begin{abstract}
\par This paper presents a new estimator of the global regularity index of a multifractional 
Brownian motion. Our estimation method is based upon a ratio statistic, which compares the realized 
global quadratic variation of a multifractional Brownian motion at two different frequencies. 
We show that a logarithmic transformation of this statistic converges in probability to the minimum of the Hurst function, which is, under weak assumptions, identical to the global regularity index of the path.    
			\end{abstract}

\setcounter{tocdepth}{4}     
\setcounter{secnumdepth}{4}  

\keywords
\newline


\section{Introduction} \label{sec1}
\setcounter{equation}{0}
\renewcommand{\theequation}{\thesection.\arabic{equation}}
\hspace{0.75cm} Fractional Brownian motion (fBm) is one of the most prominent Gaussian processes in the probabilistic and statistical literature. Popularized by Mandelbrot and van Ness \cite{Mandelbrot1968} in 1968, it found various applications in modeling stochastic phenomena in physics, biology, telecommunication and finance among many other fields. Fractional Brownian motion is characterized by its self-similarity property, the stationarity of its increments and by its ability to match any prescribed constant local regularity.
Mathematically speaking, for any $H\in (0,1)$, a fBm with Hurst index $H$, denoted by $B^H=(B^H_t)_{t \geq 0}$, is a zero mean Gaussian process with the covariance function given by
\begin{equation*} 
\E[B^H_s B^H_t] = \frac{1}{2} \left(t^{2H} + s^{2H} - |t-s|^{2H}\right).
\end{equation*} 

Various representations of fBm can be found in the existing literature; we refer to \cite{Nu,LLVH} and references therein. The Hurst parameter $H \in (0,1)$ determines the path properties of the fBm: (i) 
The process $(B^H_t)_{t \geq 0}$ is self-similar with index $H$, i.e. $(a^HB^H_t)_{t \geq 0}= 
(B^H_{at})_{t \geq 0}$ in distribution, (ii) $(B^H_t)_{t \geq 0}$ has H\"older continuous paths of any order strictly smaller than $H$, (iii) fractional Brownian motion has short memory if an only if $H \in 
(0,1/2]$. Moreover, fBm presents long range dependance if $H$ belongs to $(1/2,1)$. The statistical estimation of the Hurst parameter $H$ in the high frequency setting, i.e. the setting of mesh converging to 0 while the interval length remaining fixed,  is often performed by using 
power variation of $B^H$. Recall that a standard power variation of an auxiliary process $(Y_t)_{t
\geq 0}$ on the interval $[0,T]$ is defined by
\begin{equation*}
V(Y,p)_T^n := \sum_{i=0}^{[nT]} \left|Y_{\frac{i+1}{n}} - Y_{\frac{i}{n}} \right|^p. 
\end{equation*}
This type of approach has been investigated in numerous papers; we refer to e.g. 
\cite{GL89, IL97} among many others.  The fact that most of the properties of fBm are governed by the single parameter  $H$ restricts its application in some 
situations. In particular, its H\"older exponent remains the same along all its trajectories. This does not seem to 
be adapted to describe adequately natural terrains, for instance. In addition, long range dependence requires 
$H >1/2$, and thus imposes paths smoother than the ones of Brownian motion. Multifractional Brownian 
motion (mBm) was introduced to overcome these limitations. Several definitions of multifractional Brownian motion 
 exist. The first ones were proposed in  \cite{PL} and  \cite{ABSJDR}. A more general approach was 
introduced in \cite{StoTaq} while the most recent definition of mBm (which contains all the 
previous ones) has been given in \cite{LLVH}. The latter definition is both more flexible and retains the 
essence of this class of Gaussian processes. Recall first that a fractional Brownian field on $\R_{+}
\times(0,1)$ noted $\sifbm={(\sifbm  (t,H))}_{(t,H) \in \R_{+}\times(0,1)}$ is a Gaussian field such that, for any 
$H$, the process  ${(\sifbm  (t,H))}_{t \in \R_{+}}$ is a fBm with Hurst parameter $H$. A multifractional Brownian 
motion is simply a ``path'' traced on a fractional Brownian field. More precisely, it has been defined in 
\cite[Definition1.2.]{LLVH} as follows:
\begin{defi}
\label{def:field}
Let $h:\R_{+} \rightarrow (0,1)$ be a deterministic function and $\sifbm$ be a fractional Brownian field. A multifractional Brownian motion (mBm) with functional parameter $h$ is the Gaussian process $B^h={(B^h_t)}_{t\in \R_{+}}$ defined by $B^h_t:=\sifbm(t,h(t))$, for all $t \in \R_{+}$.
\end{defi}
Define, for any $x$ in $(0,1)$, the positive real $c_{x}$ by setting:
\begin{equation} \label{cx}
c_x :=  {\bigg(\frac{2\cos(\pi x) \Gamma(2-2x)}
{x(1-2x)} \bigg)}^{\frac{1}{2}},
\end{equation}
where $\Gamma$ denotes the standard gamma function. For any function $h:\R_{+}\rightarrow (0,1)$, it is easy to verify that the process $B^{h}:={(B^{h}_{t})}_{t\in\R_{+}}$ defined by
\begin{equation} 
\label{mfm}
B^{h}_t = \frac{1}{c_{h(t)}} \int_{\R} \frac{\exp(itx) -1}{|x|^{h(t)+1/2}} ~ W(dx),
\end{equation}
where $W$ denotes a complex Gaussian measure\footnote{See \cite{StoTaq} and \cite[Chapter $6$]{TaqSam} for more details on Gaussian complex measures.}\hspace{-0.5ex}, is a multifractional Brownian motion with functional parameter $h$.

Intuitively speaking, the  multifractional Brownian motion behaves \textit{locally} as fractional Brownian 
motion, but the functional parameter $h$ is time-varying. Moreover, it remains linked 
to local regularity of $B^{h}$, but in a less simple way than in the case of the fBm. More precisely, if we assume that $h$ belongs to the set $C^{\eta}([0,1],\R)$, for some $\eta>0$,
and is such that
\begin{equation}
 \label{hmin}
0< h_{\min}:= \min_{t \in [0,1]} h(t) \leq h_{\max}:=\max_{t \in [0,1]} h(t) < \min\{1, \eta\},
\end{equation}
then $h_{\min}$ is the regularity parameter of $B^{h}$ (see \cite[Corollaries 1,2 and Proposition 10]{ACLVL}). 
In this setting the functional parameter $h$ needs to be estimated locally in order to get a full understanding 
of the path properties of the multifractional Brownian motion $B^{h}$. Bardet and Surgailis \cite{BS13} have proposed to use a local power variation of higher order filters of increments of $B^{h}$ to estimate
the function $h$.  More specifically, they prove the law of large numbers and a central limit theorem for
the local estimator of $h$ (i) based on log-regression of the local quadratic variation, (ii) based 
on a ratio of local quadratic variations.  

In this paper we are aiming at the estimation of the parameter $h_{\min}$, which represents the
regularity (or smoothness) of the multifractional Brownian motion $B^{h}={(B^{h}_{t})}_{t\geq 0}$. For this particular statistical problem the local estimation approach investigated in 
\cite{BS13} appears to be rather inconvenient. Instead our method relies on a ratio 
statistic, which compares the \textit{global} quadratic variation at two different 
frequencies. We remark that in general it is impossible to find a global rate $a_n$ 
such that the normalized power variation $a_n V(B^{h},p)_T^n$ converges to a non-trivial 
limit. However, ratios of global
power variations can very well be useful for statistical inference. Indeed, we will show that under
appropriate conditions on the functional parameter $h$, the convergence
\begin{equation*} 
S_n(B^{h}) := \frac{\sum_{i=0}^{n-1} \left(B^{h}_{\frac{i+1}{n}} - B^{h}_{\frac{i}{n}}\right)^2}
{\sum_{i=0}^{n-2} \left(B^{h}_{\frac{i+2}{n}} - B^{h}_{\frac{i}{n}}\right)^2} \underset{n\to+\infty}{\longrightarrow} 2^{-2h_{\min}}, \qquad \text{holds in probability.}
\end{equation*}   
Then a simple log transformation gives a consistent estimator of the global regularity $h_{\min}$ of a mBm. 

The paper is structured as follows. Section 2 presents the basic distribution properties of the multifractional Brownian motion, reviews the estimation methods from \cite{BS13} and states the main 
asymptotic results of the paper. Proofs are given in Section 3.

\section{Background and main results} \label{sec2}
\setcounter{equation}{0}
\renewcommand{\theequation}{\thesection.\arabic{equation}}

In \cite{BS13} Bardet and Surgailis deal with a little bit more general processes than multifractional Brownian motions. However, in order not to overload the notations we will focus in this paper on the normalized multifractional Brownian motion (i.e.\hspace{-1.5ex} the mBm defined by \eqref{mfm}). From now on we will refer to this process as the multifractional Brownian motion and denote it 
by $B^{h}=(B^{h}_t)_{t \geq 0}$.

\subsection{Basic properties and local estimation of the functional parameter $h$}
\label{introintro}

We start with the basic properties of the mBm $B^{h}$ with functional parameter $h$.
Its covariance function is given by the expression
\begin{equation}  
\label{cov} 
R_{h}(t,s):= \E[B^{h}_t B^{h}_s] =    \frac{ c^2_{h_{t,s}}}{2c_{h(t)}c_{h(s)}}   \left( {|t|}^{2h_{t,s}} + {|s|}^{2h_{t,s}}  - {|t-s|}^{2h_{t,s}} \right),
\end{equation} 
where $h_{t,s} := \frac{h(t)+h(s)}{2}$ and $c_x$ has been defined in \eqref{cx}.
It is easy to check that $x\mapsto c_{x}$ is a $C^{\infty}((0,1))$-function.
The local behaviour of the multifractional Brownian motion is best understood via the relationship
\begin{align*}
\left(u^{-h(t)} (B^{h}_{t+us} -B^{h}_t) \right)_{s\geq 0}    \fdd \left(B^{h(t)}_s \right)_{s\geq 0}  \qquad 
\text{as } u \to 0,
\end{align*}
where $\fdd$ denotes the convergence of finite dimensional distributions. Hence, in the neighbourhood 
of any $t$ in $(0,1)$, the mBm $B^{h}$ behaves as fBm with Hurst parameter $h(t)$. This observation is essential for the local estimation of the functional parameter $h$. 
In the following we will briefly review the statistical methods of local inference investigated in  
Bardet and Surgailis \cite{BS13}, which is based on high frequency observations 
$B^{h}_0, B^{h}_{1/n}, \ldots, B^{h}_{(n-1)/n}, B^{h}_{1}$. While the original paper is investigating rather general Gaussian models whose tangent process is a fractional Brownian motion, we will specialize their asymptotic results to the framework of multifractional Brownian motion. 

Let us introduce the generalized  increments of a process $Y={(Y_{t})}_{t\geq 0}$. Consider a vector of coefficients
$a=(a_0,\ldots,a_q) \in \R^{q+1}$ and a natural number $m \geq 1$ such that 
\[
\sum_{j=0}^q  j^k a_j= 0 \quad \text{for } k=0,\ldots, m-1 \qquad \text{and} \qquad \sum_{j=0}^q  j^m a_j \not = 0.
\] 
In this case the vector $a\in \R^{q+1}$ is called a filter of order $m$. The generalised increments 
of $Y$ associated with filter $a$ at stage $i/n$ are defined as
\[
\Delta_{i,a}^n Y := \sum_{j=0}^q  a_j Y_{\frac{i+j}{n}}.
\]
Standard examples are $a^{(1)}=(-1,1)$, $\Delta_{i,a^{(1)}}^n Y=Y_{(i+1)/n} - Y_{i/n}$ (first order 
differences) and $a^{(2)}=(1,-2,1)$, $\Delta_{i,a^{(2)}}^n Y=Y_{(i+2)/n} -2Y_{(i+1)/n} +2Y_{i/n}$
(second order differences). In both cases we have that $q=m$. Now, we set 
$\psi(x,y):= (|x+y|)/(|x|+|y|)$ and set 

\[
\Lambda (H):= \E[\psi(\Delta_{0,a}^n B^{H}, \Delta_{1,a}^n B^{H})], \qquad H \in (0,1). 
\] 
The function $\Lambda$ does not depend on $n$ and  is strictly increasing on the interval $(0,1)$. For
any $\alpha \in (0,1)$, which determines the local bandwidth, the ratio type estimator of $h(t)$ is defined as 
\begin{equation}
\widehat{h}_t^{n,\alpha} :=  \Lambda^{-1} \bigg( \frac{1}{\card \{i\in \llbracket 0,n-q-1\rrbracket: |i/n-t| \leq n^{-\alpha}\}}
\sum_{i\in \llbracket 0,n-q-1\rrbracket:~|i/n-t| \leq n^{-\alpha}} \hspace{-5ex} \psi(\Delta_{i,a}^n B^{h}, \Delta_{i+1,a}^n B^{h})
\bigg).
\end{equation}
Here and throughout the paper we denote $\llbracket p,q\rrbracket:=\{p,p+1,p+2,\ldots, q\}$
for any $p,q \in \mathbb N$ with $p\leq q$.
The authors of \cite{BS13} only investigate the estimator $\widehat{h}_t^{n,\alpha}$ relative to the
filter $a=a^{(2)}$, which we assume in this subsection from now on. The consistency and asymptotic
normality of the estimator $\widehat{h}_t^{n,\alpha}$ is summarized in the following theorem. We remark that the condition for the central limit theorem crucially depends on the interplay between 
the bandwidth parameter $\alpha$ and the H\"older index $\eta$ of the function $h$.

\begin{theo} \label{th1}
(\cite[Proposition 3]{BS13}) Assume that $h$ belongs to $C^{\eta}([0,1])$ and that condition \eqref{hmin} is satisfied. 
\begin{itemize}
\item[(i)] For any $t\in (0,1)$ and $\alpha \in (0,1)$ it holds that 
\[
\widehat{h}_t^{n,\alpha} \pn h(t), \qquad \text{as } n \to \infty.
\]
\item[(ii)] When $\alpha> \max \left(\frac{1}{1+2\min(\eta,2)} , 1 - 4(\min(\eta,2) - \sup_{t \in (0,1)} h(t))\right)$ it holds that 
\[
\sqrt{2n^{1-\alpha}} \left( \widehat{h}_t^{n,\alpha} - h(t) \right) \weak \mathcal{N}(0, \tau^2) 
\qquad \text{as } n \to \infty,
\]
where the asymptotic variance $\tau^2$ is defined in \cite[Eq. (2.17)]{BS13}.
\end{itemize}
\end{theo} 
The paper \cite{BS13} contains the asymptotic theory for a variety of other local estimators of $h(t)$.
We dispense with the detailed exposition of these estimators, since only  $\widehat{h}_t^{n,\alpha}$ 
is somewhat related to our estimation method.  

\begin{rem} \label{rem1} \rm
Nowadays, it is a standard procedure to consider higher order filters for Gaussian processes to obtain 
a central limit theorem for the whole range of Hurst parameters. Let us shortly recall some classical 
asymptotic results, which are usually referred to as Breuer-Major central limit theorems. We consider
the scaled power variation of a fractional Brownian motion $B^H$ with Hurst parameter $H \in (0,1)$ based 
on first order filter $a^{(1)}$ and second order filter $a^{(2)}$:
\[
V(B^H,p; a^{(1)})^n := n^{-1+pH}\sum_{i=0}^{n-1} |\Delta_{i,a^{(1)}}^n B^H|^p
\hspace{2.5ex} \text{and} \hspace{2.5ex} 
V(B^H,p; a^{(2)})^n := n^{-1+pH}\sum_{i=0}^{n-2} |\Delta_{i,a^{(2)}}^n B^H|^p.
\]
It is well known that, after an appropriate normalization, 
the statistic $V(B^H,p; a^{(1)})^n$ exhibits asymptotic normality for $H \in (0,3/4]$, while it converges 
to the Rosenblatt distribution for $H \in (3/4,1)$. On the other hand, the statistic $V(B^H,p; a^{(2)})^n$ 
exhibits asymptotic normality for all $H \in (0,1)$. We refer to \cite{BM83,T79} for a detailed exposition. 
\qed
\end{rem}

\subsection{Estimation of the global regularity parameter $h_{\min}$}
In this section we will construct a consistent estimator of the global regularity parameter $h_{\min}$,
which has been defined at  \eqref{hmin}. Our first condition is on the set $h^{-1}(\{h_{\min}\})$, which is necessarily compact since $h$ belongs to $C^{\eta}([0,1])$. We assume that this set has the following form
\begin{align} \label{mh}
\ds {\cM}_{h}:=h^{-1}(\{h_{\min}\})= \left(\overset{q}{\underset{i=1}{\bigcup}} [a_{i},b_{i}] 
\right) \bigcup \ \left(\overset{m}{\underset{j=1}{\bigcup}} \{x_{j}\}\right), \qquad
(q,m) \in \N^2 \setminus (0,0),
\end{align} 
where $\N=\{0,1,2,\ldots\}$ and the intervals $ [a_{i},b_{i}] $ are disjoint and such that 
 none of the $x_{j}$'s belongs to  $\overset{q}{\underset{i=1}{\bigcup}} 
[a_{i},b_{i}]$. Depending on whether $q\geq 1$ or $q=0$, we will need an additional assumption.
Below, we denote by $h^{(p)}_l(x)$ (resp. $h^{(p)}_r(x)$) the $p$th left (resp. right) derivative
of $h$ at point $x$.  \\ \\  \\
($\mathscr{A}$)
 There exist  positive integers $p_j$ such that function $h$  is $p_j$ times continuously left 
 and right differentiable at point $x_j$ for $j=1,\ldots,m$ such that 
 \[
p_j = \min\{p: ~ h^{(p)}_l(x_j)\not= 0\} = \min\{p: ~ h^{(p)}_r(x_j) \not= 0\}. 
 \]
 \newline
We remark that since $h$ reaches its minimum at points $x_j$, we necessarily have 
that $h^{(p_j)}_r(x_j) > 0$ and that $h^{(p_j)}_l(x_j)>0$ 
if $p_j$ is even and $h^{(p_j)}_l(x_j)<0$ if $p$ is odd.
Now, we proceed with the construction of the consistent estimator of  
the global regularity parameter $h_{\min}$ based on high frequency observations $B^{h}_0, B^{h}_{1/n}, \ldots, 
B^{h}_{(n-1)/n}, B^{h}_{1}$. 
First of all, let us remark that considering the estimator $\min_{t \in [0,1]} 
\widehat{h}_t^{n,\alpha}$, where $\widehat{h}_t^{n,\alpha}$ has been introduced 
in the previous section, is not a trivial matter since the functional version of Theorem \ref{th1} is not 
available. Instead our statistics relies on the global quadratic variation rather than local estimates. 

For the mBm $B^h={(B^h_{t})}_{t\in [0,1]}$, we introduce the notations
\begin{align} \label{vnsn}
V(B^h;k)^n:= \sum_{i=0}^{n-k} \left(B^h_{\frac{i+k}{n}} - B^h_{\frac{i}{n}}\right)^2, \qquad S_n(B^h):= \frac{V(B^h;1)^n}{V(B^h;2)^n}.
\end{align} 
Our first result determines the limit of $\E[V(B^{h};1)^n] / \E[V(B^{h};2)^n]$.

\begin{prop} \label{firstresult}
Let $h:[0,1]\rightarrow (0,1)$ be a deterministic $C^{\eta}([0,1])$-function satisfying \eqref{hmin} 
 and such that 
the set $ {\cM}_{h}$ has the form \eqref{mh}. If $q=0$ we also assume that condition ($\mathscr{A}$) holds.
Define 
\[
\ds \uh:=\frac{\E[V(B^{h};1)^n]}{\E[V(B^{h};2)^n]}.
\]
Then it holds that 
\begin{align} \label{first}
\underset{n\to+\infty}{\lim}\  \uh = {\left(\frac{1}{2}\right)}^{2h_{\min}}.
\end{align}
\end{prop}

The convergence result of Proposition \ref{firstresult} is rather intuitive when $q\geq 1$, which means
that the minimum of the function $h$ is reached on a set of positive Lebesgue measure. In this setting 
it is quite obvious that the statistic $V(B^{h};k)^n$ is dominated by squared increments 
$(B^{h}_{(i+k)/n} - B^{h}_{i/n})^2$ for $i/n \in \cup_{i=1}^q [a_{i},b_{i}] $. Thus, the estimation problem is similar to the estimation of the Hurst parameter of a fractional Brownian motion 
$(B_t^{h_{\min}})_{t\in \cup_{i=1}^q [a_{i},b_{i}]}$ with Hurst parameter  $h_{min}$, 
for which the convergence at 
\eqref{first} is well known. When $q=0$, and hence $\text{Leb}( {\cM}_{h})=0$, 
the proof of Proposition \ref{firstresult} becomes much more delicate.

\begin{rem} \rm
Assume for illustration purpose 
that $q=0$, $m=1$, $x:=x_1$ and $p:=p_1$. Condition ($\mathscr{A}$) is crucial 
to determine the precise asymptotic expansion of the quantity $\E[V(B^{h};k)^n]$. The lower
and upper bounds in \eqref{oifjefoijefijoreijr1dede} and \eqref{oifjefoijefijoreijr2dede} in the proof show that 
\[
\E[V(B^{h};k)^n]= \text{O}\left( \frac{n^{1-2h_{\min}}}{(\ln n)^{1/p}}\right)
\quad \text{as } n\to + \infty, \qquad \text{for } k=1,2.
\]   
The condition $ \min\{p: ~ h^{(p)}_l(x)\not= 0\} = \min\{p: ~ h^{(p)}_r(x) \not= 0\}$ of assumption 
 ($\mathscr{A}$) is not essential for the proofs. For instance, when 
 $ \min\{p: ~ h^{(p)}_l(x)\not= 0\} > \min\{p: ~ h^{(p)}_r(x) \not= 0\}$ the expectation $\E[V(B^{h};k)^n]$
 would be dominated by the terms in the small neighbourhood on the right hand side of $x$
 and the statement of Proposition \ref{firstresult} can be proved in the same manner.  \qed
\end{rem}
Our main result shows that the statistic $S_n(B^{h})$ and the quantity $\ds \uh$ are asymptotically equivalent
in probability.

\begin{theo}
 \label{th2}
 Assume that $h\in C^2([0,1])$ and 
the set $ {\cM}_{h}$ has the form \eqref{mh}. If $q=0$ we also assume that condition ($\mathscr{A}$) holds.
Then we have the following result:
 \begin{equation} \label{asyequi}
S_n(B^{h}) \pn {\left(\frac{1}{2}\right)}^{2h_{\min}}.
\end{equation}
In particular, the following convergence holds:  
\begin{equation} \label{hminhat}
\widehat{h}_{\min} := -\frac{\ln (S_n(B^{h}))}{2\ln (2)} \pn h_{\min}.
\end{equation}
\end{theo}

%
%

The asymptotic result of Theorem \ref{th2} can be extended to more general Gaussian processes than 
the mere multifractional Brownian motion. As it has been discussed in  \cite{BS13}, when a Gaussian process possesses a tangent process $B^{h(t)}$ at time $t$, we may expect Theorem \ref{th2} to hold under certain assumptions on its covariance kernel. We refer to assumptions $(A)_{\kappa}$
and $(B)_{\alpha}$ therein for more details on sufficient conditions.   

The rate of convergence, or a weak limit theorem, associated with the consistency result at \eqref{hminhat} is a more delicate issue. In the setting $q=0$, which implies that $\text{Leb}( {\cM}_{h})=0$, the bias associated with the convergence in  \eqref{first} may very well dominate the variance 
of the estimator. A careful inspection of our proof, and more specifically of statement
\eqref{boundmu} and Remark \ref{ofjerifrofjerfioerj}, implies that the bias in Proposition \ref{firstresult} has a logarithmic rate.   Thus, weak limit theorems for the estimator  $\widehat{h}_{\min}$ are out of reach in this framework. When $q \geq 1$ and hence $\text{Leb}( {\cM}_{h})>0$, one may hope to find better rates 
of convergence for the estimator $\widehat{h}_{\min}$. However, we dispense with the exact exposition of this statistical problem.

\section{Proofs} \label{sec3}
\setcounter{equation}{0}
\renewcommand{\theequation}{\thesection.\arabic{equation}}

Throughout this section we denote all positive constants by $C$, or $C_p$ if they depend on an external parameter $p$, although they may change from line to line.

\subsection{Proof of Proposition \ref{firstresult}}
For $k=1,2$ we introduce the notation 
\begin{align} \label{vkn}
V_n^{(k)}:= \sum_{i=0}^{n-k} \left(\frac{k}{n} \right)^{2h(i/n)},
\end{align}
which serves as the first order approximation of the quantity $\E[V(B^{h};k)^n]$.
Applying \cite[Lemma 1 p.$13$]{BS10} we conclude that 
\begin{equation} 
\label{approx1}
\left| \E[V(B^{h};k)^n] -  V_n^{(k)} \right| \leq C \frac{\ln n}{n^{\eta\wedge 1}}  
\sum_{i=0}^{n-k} \left(\frac{i}{n} \right)^{2h(k/n)} \leq C
\frac{\ln n}{n^{2h_{\min} - 1 +\eta\wedge 1}} 
\end{equation}
for any $(n,k) \in \N \times \{1,2\}$. We have the inequality 
\begin{align} 
\left|  \uh  - \left(\frac{1}{2}\right)^{2h_{\min}}  \right|
&\leq \frac{|\E[V(B^{h};1)^n] -  V_n^{(1)} | +  |\E[V(B^{h};2)^n] -  V_n^{(2)} |}{V_n^{(2)}}
+ \left| \frac{V_n^{(1)} }{V_n^{(2)}} - \left(\frac{1}{2}\right)^{2h_{\min}}   \right| \nonumber \\[1.5 ex]
\label{def1} &=:\mu^{(1)}_{n} + \mu^{(2)}_{n}. 
\end{align}
We first show that $\mu^{(1)}_{n} \to 0$ as $n\to \infty$. When $h_{\min}= h_{\max}$ we trivially 
have $\mu^{(1)}_{n} =0$. If $h_{\min}< h_{\max}$,  we fix $\epsilon \in (0,h_{\max} - h_{\min})$. By 
$\text{Leb}(A)$ we denote the Lebesgue measure of any measurable set $A$. We have
that 
$$\text{Leb}\left(h^{-1}([h_{\min}, h_{\min} + \epsilon])\right)>0.$$
Thus, there exists $n_0 \in \N$ such that for all $n\geq n_0$ it holds that 
\[
\text{Card}\{i\in \llbracket 0,n-k\rrbracket;~ h(i/n) \in [h_{\min}, h_{\min} + \epsilon] \} \geq n\ \text{Leb}\left(h^{-1}([h_{\min}, h_{\min} + \epsilon])\right)/2. 
\] 
This implies that 
\begin{align*}
V_n^{(2)} \geq  \sum_{i\in \llbracket 0,n-k\rrbracket;  ~h(i/n) \in [h_{\min}, h_{\min} + \epsilon] }\left(\frac{2}{n} \right)^{2h(i/n)}
\geq C n^{1 - 2(h_{\min} + \epsilon)}.
\end{align*}
Hence, applying Inequality \eqref{approx1}, we conclude that:
\[
\mu^{(1)}_{n} \leq C \ln n\cdot n^{2\epsilon  -\eta\wedge 1},
\]
which proves that $\mu^{(1)}_{n} \underset{n\to+\infty}{\to } 0$, for any $\epsilon$ small enough. \qed
%

\subsubsection{Convergence of 
$\mu^{(2)}_{n}$ in the case $q\geq 1$}
We first prove that $\mu^{(2)}_{n} \to 0$ in the case $q\geq 1$. Assume again that 
$h_{\min}< h_{\max}$. First, we observe the lower bound
\begin{align} 
V_n^{(k)} &\geq  \sum_{l=1}^q~ \sum_{i\in {\tiny \llbracket 0,n-k\rrbracket}; ~ i/n \in [a_l,b_l] }\left(\frac{k}{n} \right)^{2h(i/n)}
= \left( \frac{k}{n}\right)^{2h_{\min}}   \sum_{l=1}^q 
\card\{i\in \llbracket 0,n-k\rrbracket; ~ i/n \in [a_l,b_l] \} \nonumber \\[1.5 ex]
\label{lowercase1} & \geq
 n \left( \frac{k}{n}\right)^{2h_{\min}} \sum_{l=1}^q \left(b_l-a_l-\tfrac{2}{n} \right). 
\end{align}

For the upper bound we fix $0<\epsilon< h_{\max} - h_{\min}$ and consider the decomposition 
\begin{align*}
V_n^{(k)}= \sum_{i\in {\llbracket 0,n-k\rrbracket}; ~h(i/n) \in [h_{\min}, h_{\min} + \epsilon] }\left(\frac{k}{n} \right)^{2h(i/n)}
+ \sum_{i\in {\llbracket 0,n-k\rrbracket}; ~h(i/n) \not \in [h_{\min}, h_{\min} + \epsilon] }\left(\frac{k}{n} \right)^{2h(i/n)}.
\end{align*}
Setting $\lambda_n(\epsilon) := n^{-1}\text{card}\{i\in \llbracket 0,n-k\rrbracket; ~h(i/n) \in [h_{\min}, h_{\min} + \epsilon]\}$, 
we deduce the assertions
\begin{align*}
\lambda_n(\epsilon) &\to \text{Leb}\left(h^{-1}([h_{\min}, h_{\min} + \epsilon])\right) \quad \text{as }
n \to \infty,  \\[1.5 ex]
\text{Leb}\left(h^{-1}([h_{\min}, h_{\min} + \epsilon])\right) &\to 
\text{Leb}\left(h^{-1}(\{h_{\min}\})\right) =  \sum_{l=1}^q (b_l-a_l)>0 \quad \text{as } \epsilon \to 0. 
\end{align*}
Now, we conclude that
\begin{align} \label{uppercase1}
V_n^{(k)} \leq n \lambda_n(\epsilon) \left(\frac{k}{n} \right)^{2h_{\min}} 
+ n (1-\lambda_n(\epsilon)) \left(\frac{k}{n} \right)^{2(h_{\min}+\epsilon)}. 
\end{align}
Throughout the proofs we write $\underline{\lim}$ for $\lim\inf$ and $\overline{\lim}$ for $\lim\sup$.
Applying inequalities \eqref{lowercase1} and \eqref{uppercase1}, we obtain that
\begin{multline*}
\lin \frac{n \sum_{l=1}^q (b_l-a_l-\frac{2}{n})} {n \lambda_n(\epsilon) 
+ n (1-\lambda_n(\epsilon)) \left(\frac{2}{n} \right)^{2\epsilon}}
%
  \\ 
\leq 
 \lin 2^{2h_{\min}} \frac{V_n^{(1)} }{V_n^{(2)}} \leq 
\lsp 2^{2h_{\min}} \frac{V_n^{(1)} }{V_n^{(2)}} \leq
 \\
 \lsp \frac {n \lambda_n(\epsilon) 
+ n (1-\lambda_n(\epsilon)) \left(\frac{1}{n} \right)^{2\epsilon}}
{n \sum_{l=1}^q (b_l-a_l-\frac{2}{n})}.
\end{multline*}

Hence, we deduce that  
\begin{align*}
&
 \frac{2^{-2h_{\min}}\text{Leb}\left(h^{-1}(\{h_{\min}\})\right)} {\text{Leb}\left(h^{-1}([h_{\min}, h_{\min} + \epsilon])\right)} 
\leq \lin \frac{V_n^{(1)} }{V_n^{(2)}} \leq 
\lsp \frac{V_n^{(1)} }{V_n^{(2)}} 
&
 \leq \frac{2^{-2h_{\min}}\text{Leb}\left(h^{-1}([h_{\min}, h_{\min} + \epsilon])\right)} 
{\text{Leb}\left(h^{-1}(\{h_{\min}\})\right)}. 
\end{align*}
By letting $\epsilon$ tend to $0$, we readily deduce taht $\mu^{(2)}_{n} \to 0$ as $n\rightarrow +\infty$.

\subsubsection{Convergence of 
$\mu^{(2)}_{n}$ in the case $q= 0$}
Without loss of generality we assume that $m=1$ and ${\cM}_{h}=h^{-1}(\{h_{\min}\})=\{x\}$
 with $x\in (0,1)$. Recall that
in this setting we assume condition ($\mathscr{A}$) with $p:=p_1$. We let $\gamma$ be a
positive number such that 
$\gamma<2^{-1}\min\{{|h^{(p)}_{l}(x)|,h^{(p)}_{r}(x)}\}$. Now, there exists a $\epsilon
=\epsilon(\gamma)>0$ with $\epsilon< \min\{x,1-x, \gamma\}$ such that: 
\begin{align}
 \label{deriv1}
&\hspace{-17ex}\forall y>x \text{ with } 0<y-x<\epsilon: \notag \\[1.5 ex]
&\hspace{-17ex}h_{\min} +\frac{1}{p!} (y-x)^{p}\ (h^{(p)}_{r}(x) - \gamma) \leq h(y) \leq h_{\min} + 
\frac{1}{p!} (y-x)^{p} \ (h^{(p)}_{r}(x) + \gamma),
\end{align}
\vspace{-3ex}
\begin{align}
\label{deriv2}
&\forall y<x \text{ with } 0<x-y<\epsilon: \notag \\[1.5 ex]
&h_{\min} + \frac{1}{p!}  (y-x)^{p}\ (h^{(p)}_{l}(x) -(-1)^{p} \gamma) \leq h(y) \leq h_{\min} + 
\frac{1}{p!}  (y-x)^{p}\ (h^{(p)}_{l}(x) +(-1)^{p} \gamma).
\end{align}
We proceed with the derivation of upper and lower bounds for the quantity $\mu^{(2)}_{n}$. 
We start
with the decomposition $V_n^{(k)} = \Gamma^{(1)}_{n,k}(\gamma, \epsilon) + \Gamma^{(2)}_{n,k}(\gamma, \epsilon) +\Gamma^{(3)}_{n,k}(\gamma, \epsilon)$ where
\begin{align*}
\Gamma^{(1)}_{n,k}(\gamma, \epsilon)&:= 
\sum_{i\in \llbracket 0,n-k\rrbracket;  ~i/n \in [x,x+\epsilon] }\left(\frac{k}{n} \right)^{2h(i/n)};\\
\Gamma^{(2)}_{n,k}(\gamma, \epsilon) &:= 
\sum_{i\in \llbracket 0,n-k\rrbracket;  ~i/n \in [x-\epsilon, x) }\left(\frac{k}{n} \right)^{2h(i/n)}; \\
\ds \Gamma^{(3)}_{n,k}(\gamma, \epsilon) &:= 
\sum_{i\in \llbracket 0,n-k\rrbracket; ~i/n \in [x-\epsilon, x+\epsilon]^c }\left(\frac{k}{n} \right)^{2h(i/n)}.
\end{align*} 
It is clear that $\Gamma^{(3)}_{n,k}(\gamma, \epsilon) \leq n (k/n)
^{2 h(y_{\varepsilon})}$, where we have set 
$$y_{\epsilon}:=\text{argmin}
\{h(u): \ u\in(x-\epsilon,x+\epsilon)^{c}\cap[0,1]  \}.$$
For the other two quantities, we deduce that $\underline{\Gamma}^{(r)}_{n,k}(\gamma, \epsilon) \leq \Gamma^{(r)}_{n,k}(\gamma, \epsilon) \leq 
\overline{\Gamma}^{(r)}_{n,k}(\gamma, \epsilon)$ with 
\begin{align*}
\underline{\Gamma}^{(1)}_{n,k}(\gamma, \epsilon) &:=\left(\frac{k}{n}\right)^
{2h_{\min}} \sum_{i\in \llbracket 0,n-k\rrbracket: ~i/n \in [x,x+\epsilon] }\left(\frac{k}{n} \right)^{2
(p!)^{-1} (i/n-x)^{p}(h^{(p)}_{r}(x) + \gamma)}, \\[1.5 ex]
\underline{\Gamma}^{(2)}_{n,k}(\gamma, \epsilon) &:=\left(\frac{k}{n}\right)^
{2h_{\min}} \sum_{i\in \llbracket 0,n-k\rrbracket: ~i/n \in [x-\epsilon,x) }\left(\frac{k}{n} \right)^{2
(p!)^{-1} (i/n-x)^{p}(h^{(p)}_{l}(x) + (-1)^{p} \gamma)}
\end{align*}
 and $\overline{\Gamma}^{(1)}_{n,k}(\gamma, \epsilon) :=\underline{\Gamma}^{(1)}_{n,k}(-\gamma, \epsilon) $ and 
$\overline{\Gamma}^{(2)}_{n,k}(\gamma, \epsilon):=\underline{\Gamma}^{(2)}_{n,k}(-\gamma, \epsilon)$. 
Using \eqref{deriv1} and \eqref{deriv2}, it is easy to see that, for every $(k,n)\in \{1,2\}
\times\bN$:
\begin{equation} 
\label{boundmu} 
\underline{\mu}^{(2)}_{n}(\gamma, \epsilon)
\leq \frac{V_n^{(1)} }{V_n^{(2)}} \leq
\overline{\mu}^{(2)}_{n}(\gamma, \epsilon),
\end{equation} 
with
\begin{align*} 
&\underline{\mu}^{(2)}_{n}(\gamma, \epsilon) :=\frac{\underline{\Gamma}^{(1)}_{n,1}(\gamma, \epsilon) 
+ \underline{\Gamma}^{(2)}_{n,1}(\gamma, \epsilon) }
{\overline{\Gamma}^{(1)}_{n,2}(\gamma, \epsilon) + 
\overline{\Gamma}^{(2)}_{n,2}(\gamma, \epsilon) 
+ n (2/n)
^{2h(y_{\varepsilon})}}, \nonumber  \\
&\overline{\mu}^{(2)}_{n}(\gamma, \epsilon) :=
\frac{\overline{\Gamma}^{(1)}_{n,1}(\gamma, \epsilon) + 
\overline{\Gamma}^{(2)}_{n,1}(\gamma, \epsilon) 
+ n (1/n)
^{2 h(y_{\varepsilon})}}
{\underline{\Gamma}^{(1)}_{n,2}(\gamma, \epsilon) 
+ \underline{\Gamma}^{(2)}_{n,2}(\gamma, \epsilon) }.  \nonumber
\end{align*} 

From \eqref{boundmu} we obtain that 
\begin{equation}
\label{dezzezedze}
0\leq 2^{2h_{\min}}  \mu^{(2)}_{n} \leq \left |2^{2h_{\min}} \frac{V_n^{(1)} }{V_n^{(2)}} -1
\right| \leq U_{n}(\gamma, \epsilon)+U_{n}(-\gamma, \epsilon),
\end{equation}
where 
\begin{align}
\label{zodzeoijzoecijs1}
U_{n}(\gamma, \epsilon)&:= |\underline{\Delta}_{n,2}(\gamma, \epsilon)|^{-1} \bigg({|2^{2h_{\min}}\overline{\Delta}_{n,1}(\gamma, \epsilon)-\underline{\Delta}_{n,2}(\gamma, \epsilon)|}  + 2 {n}^{1-2h(y_{\epsilon})}\bigg), \\[1.5 ex]
\underline{\Delta}_{n,k}(\gamma, \epsilon)&:=\underline{\Gamma}^{(1)}_{n,k}(\gamma, \epsilon)  + \underline{\Gamma}^{(2)}_{n,k}(\gamma, \epsilon), \qquad 
\overline{\Delta}_{n,k}(\gamma, \epsilon):=\underline{\Delta}_{n,k}(-\gamma, \epsilon). 
\end{align}
In view of \eqref{dezzezedze} it is  sufficient to show that $\underset{\gamma\to 0}{\overline{\lim}}\ \lsp U_{n}(\gamma, \epsilon)=0$. 
Define 
\[
d_{\gamma}:=2(p!)^{-1}(h^{(p)}_{r}(x)+\gamma) \qquad \text{and} \qquad 
d'_{\gamma}:=2(p!)^{-1}(h^{(p)}_{l}(x)+ (-1)^{p}\gamma).
\]
For any $(a,b)$ in $\R_{+}\times(\R \setminus\{0\})$, we also set 
\begin{align}
\label{edpokezpodk1}
S_{n,k}(a,\epsilon)&:=  \sum_{i\in \llbracket 0,n-k\rrbracket: ~i/n \in [x,x+\epsilon] }\left(\frac{k}{n} \right)^{a (i/n-x)^{p}}, \\
\label{edpokezpodk2}
T_{n,k}(b,\epsilon)&:= \sum_{i\in \llbracket 0,n-k\rrbracket: ~i/n \in [x-\epsilon,x) }\left(\frac{k}{n} \right)^{b (i/n-x)^{p}}.
\end{align}
We deduce the identities $\underline{\Gamma}^{(1)}_{n,k}(\gamma, \epsilon)= (k/n)^{2h_{\min}}S_{n,k}(d_{\gamma},
\epsilon)$ and $\underline{\Gamma}^{(2)}_{n,k}(\gamma, \epsilon)= (k/n)^{2h_{\min}}T_{n,k}(d'_{\gamma},
\epsilon)$. Note moreover that $d'_{\gamma}>0$ when $p$ is even and $d'_{\gamma}<0$ when $p$ is odd.
We therefore assume from now on that $b>0$ when $p$ is even and that $b<0$ when $p$ is odd. 
For any $\eta \in \R\setminus \{0\}$, we define 
$$f^{(\eta)}_{n,k}(u):=\left(\frac{k}{n} \right)^{\eta(u-x)^{p}}.$$ 
Since $i\mapsto f^{(a)}_{n,k}(i/n)$ is decreasing on $\llbracket [nx]+1, [n(x+\epsilon)] \rrbracket$ while
$i\mapsto f^{(b)}_{n,k}(i/n)$ is increasing if $p$ even (resp. decreasing if $p$ odd) on $\llbracket [n(x-\epsilon)]+1,[nx]  \rrbracket$, one can use an integral test for convergence, which provides us with the following upper bounds
\begin{align}
\label{A1}
\frac{n\int^{\beta_{n}(a)}_{\alpha_{n}(a)} y^{1/p-1} e^{-y} \ dy }{p (a\ln(n/k))^{1/p}}
&\leq 
S_{n,k}(a,\epsilon)
\leq
\frac{n\int^{\mu_{n}(a)}_{\tau_{n}(a)} y^{1/p-1} e^{-y} \ dy }{p (a\ln(n/k))^{1/p}}, \\[1.5 ex]
\label{A2}
\frac{n\left(\int^{\beta'_{n}(b)}_{\alpha'_{n}(b)} y^{1/p-1} e^{-y} \ dy-\rho^{(b)}_{n,k}(\epsilon)
\right)}{p ((-1)^{p} b \ln(n/k))^{1/p}}
&\leq 
T_{n,k}(b,\epsilon)
\leq
\frac{n\left(\int^{\mu'_{n}(b)}_{\tau'_{n}(b)} y^{1/p-1} e^{-y} \ dy-\rho^{(b)}_{n,k}(\epsilon)
\right)}{p ((-1)^{p} b\ln(n/k))^{1/p}}.
\end{align}
Here we use the notation 
\begin{align*}
\alpha_{n}(a)&:= a \ln(n/k) \bigg(\frac{[nx]+1}{n}-x\bigg)^{p}, \quad \beta_{n}(a):= a \ln(n/k) \bigg(\frac{[n(x+\epsilon)]+1}{n}-x\bigg)^{p},\\
\tau_{n}(a)&:= a \ln(n/k) \bigg(\frac{[nx]}{n}-x\bigg)^{p} , \quad \mu_{n}(a):= a \ln(n/k) \bigg(\frac{[n(x+\epsilon)]}{n}-x\bigg)^{p}
\end{align*}
and $\rho^{(b)}_{n,k}(\epsilon):=f^{(b)}_{n,k}(\frac{[n(x-\epsilon)]+1}{n})+f^{(b)}_{n,k}(\frac{[nx]}{n})$.
Furthermore,
\begin{align*}
(\alpha'_{n}(b), \beta'_{n}(b), \tau'_{n}(b),\mu'_{n}(b))&:=(z^{(1)}_{n}(b),z^{(2)}_{n}(b),z^{(3)}_{n}(b),z^{(4)}_{n}(b)) \ \text{if $p$ is even},\\
(\alpha'_{n}(b), \beta'_{n}(b), \tau'_{n}(b),\mu'_{n}(b))&:=(z^{(3)}_{n}(b),z^{(4)}_{n}(b),z^{(1)}_{n}(b),z^{(2)}_{n}(b)) \ \text{if $p$ is odd,}
\end{align*}
where we have set
\begin{align*}
z^{(1)}_{n}(b)&:= b \ln(n/k) \bigg(\frac{[nx]-2}{n}-x\bigg)^{p}, \quad z^{(2)}_{n}(b):= b \ln(n/k) \bigg(\frac{[n(x-
\epsilon)]+1}{n}-x\bigg)^{p},\\
z^{(3)}_{n}(b)&:= b \ln(n/k) \bigg(\frac{[nx]-1}{n}-x\bigg)^{p}, \quad z^{(4)}_{n}(b):= b \ln(n/k) \bigg(\frac{[n(x-
\epsilon)]+2}{n}-x\bigg)^{p}.
\end{align*}
In view of the inequalities  \eqref{A1} and \eqref{A2}, as well as identities \eqref{edpokezpodk1} and \eqref{edpokezpodk2}, we then deduce that 
\begin{align}
\label{oifjefoijefijoreijr1dede}
&\frac{n^{1-2h_{\min}} k^{2h_{\min}} u_{n,k,p}(d_{\gamma})}{(\ln (n/k))^{1/p}}\cdot \bigg(\frac{1}{d_{\gamma}}\bigg)
&\leq& 
&\underline{\Gamma}^{(1)}_{n,k}(\gamma, \epsilon) 
&\leq&
&\frac{n^{1-2h_{\min}} k^{2h_{\min}}v_{n,k,p}(d_{\gamma})}{(\ln (n/k))^{1/p}}\cdot \bigg(\frac{1}{d_{\gamma}}\bigg) , \\
\label{oifjefoijefijoreijr2dede}
&\frac{n^{1-2h_{\min}} k^{2h_{\min}}u'_{n,k,p}(d'_{\gamma})}{(\ln (n/k))^{1/p}}\cdot \bigg(\frac{1}{|d'_{\gamma}|}\bigg)
&\leq&
&\underline{\Gamma}^{(2)}_{n,k}(\gamma, \epsilon) 
&\leq&
&\frac{n^{1-2h_{\min}} k^{2h_{\min}}v'_{n,k,p}(d'_{\gamma})}{(\ln (n/k))^{1/p}}\cdot \bigg(\frac{1}{|d'_{\gamma}|}\bigg). 
\end{align}
Here we have used the notation
\begin{align*}
u_{n,k,p}(a)&:= \frac{1}{p} \int^{{\beta_{n}(a)}}_{\alpha_{n}(a)} \ y^{1/p-1} e^{-y}\ dy, \hspace{4ex} v_{n,k,p}(a):= \frac{1}{p} \int^{{\mu_{n}(a)}}_{\tau_{n}(a)} \ y^{1/p-1} e^{-y}\ dy,\\
u'_{n,k,p}(b)&:= \frac{1}{p} \int^{{\beta'_{n}(b)}}_{\alpha'_{n}(b)} \ y^{1/p-1} e^{-y}\ dy-\frac{{\big((-1)^{p} b \ln(n/k)\big)}^{1/p} \rho^{(b)}_{n,k}(\epsilon)}{pn}, \\
v'_{n,k,p}(b)&:= \frac{1}{p} \int^{{\mu'_{n}(b)}}_{\tau'_{n}(b)} \ y^{1/p-1} e^{-y}\ dy-\frac{{\big((-1)^{p} b \ln(n/k)\big)}^{1/p} \rho^{(b)}_{n,k}(\epsilon)}{pn}.
\end{align*}
Since $\underline{\Gamma}^{(r)}_{n,k}(\gamma, \epsilon) = \overline{\Gamma}^{(r)}_{n,k}(-\gamma, \epsilon)$, \eqref{oifjefoijefijoreijr1dede} and \eqref{oifjefoijefijoreijr2dede} also provide us with upper and lower bounds for $\overline{\Gamma}^{(r)}_{n,k}(\gamma, \epsilon)$. Finally, we obtain the following lower and upper bounds
\begin{align}
\label{zdozepodkezpodkzepodkezp}
\frac{n^{1-2h_{\min}} k^{2h_{\min}}}{(\ln (n/k))^{1/p}}\cdot \Lambda_{n,k}(\gamma, \epsilon)
&\leq
\underline{\Delta}_{n,k}(\gamma, \epsilon)
\leq
\frac{n^{1-2h_{\min}} k^{2h_{\min}}}{(\ln (n/k))^{1/p}} \Lambda'_{n,k}(\gamma, \epsilon),\\
\frac{n^{1-2h_{\min}} k^{2h_{\min}}}{(\ln (n/k))^{1/p}}\cdot \Lambda_{n,k}(-\gamma, \epsilon)
&\leq
\overline{\Delta}_{n,k}(\gamma, \epsilon)
\leq
\frac{n^{1-2h_{\min}} k^{2h_{\min}}}{(\ln (n/k))^{1/p}} \Lambda'_{n,k}(-\gamma, \epsilon),
\end{align}
where
\begin{align*}
\Lambda_{n,k}(\gamma, \epsilon)&:= \frac{1}{d_{\gamma}}\cdot u_{n,k,p}(d_{\gamma}) + \frac{1}{|d'_{\gamma}|} \cdot u'_{n,k,p}(d'_{\gamma}),\\
\Lambda'_{n,k}(\gamma, \epsilon)&:= \frac{1}{d_{\gamma}}\cdot v_{n,k,p}(d_{\gamma}) + \frac{1}{|d'_{\gamma}|} \cdot v'_{n,k,p}(d'_{\gamma}).
\end{align*}
Denote $c_{p}:=\int^{+\infty}_{0} y^{1/p-1} e^{-y}\ dy$.  Recalling the definition of the constants 
$d_{\gamma}$ and $d'_{\gamma}$, a straightforward computation shows that, for any $(k,k')\in {\{1,2\}}^{2}$ with $k\neq k'$: 
\begin{align}
\label{zoiejozedijezoidjzeodijzeodijzzz1}
&\underset{n\to+\infty}{\lim} \Lambda_{n,k}(\gamma, \epsilon) =  \underset{n\to+\infty}{\lim} \Lambda'_{n,k}(\gamma, \epsilon) = \frac{c_{p}}{p}(1/d_{\gamma} + 1/|d'_{\gamma}|), \\[1.5 ex]
\label{zoiejozedijezoidjzeodijzeodijzzz2}
&  \lsp |\Lambda'_{n,k'}(\gamma, \epsilon)- \Lambda_{n,k}(-\gamma, \epsilon)| \leq C\ (2 |\gamma| + |1/d_{-\gamma}-1/d_{\gamma} + 1/|d'_{-\gamma}| - 1/|d'_{\gamma}||)\leq C |\gamma|.
\end{align}
Starting from \eqref{zdozepodkezpodkzepodkezp}, and using \eqref{zoiejozedijezoidjzeodijzeodijzzz1} and \eqref{zoiejozedijezoidjzeodijzeodijzzz2}, we see that there exists a positive integer $n_{0}$ and $C>0$  such that for all $n\geq n_{0}$
\begin{equation}
\label{dzeoijdzeoidjzeoidjzeoidjzeodjoizejdoziejdoeizjdeiz}
|\underline{\Delta}_{n,k}(\gamma, \epsilon)|^{-1}
\leq C
\frac{(\ln (n/k))^{1/p}}{n^{1-2h_{\min}} k^{2h_{\min}}}. 
\end{equation}
Finally, inequalities \eqref{zoiejozedijezoidjzeodijzeodijzzz1}, \eqref{zoiejozedijezoidjzeodijzeodijzzz2} and \eqref{dzeoijdzeoidjzeoidjzeoidjzeodjoizejdoziejdoeizjdeiz} imply that there exists a positive integer $N$ such that for all $n\geq N$:
\begin{equation*}
U_{n}(\gamma, \epsilon)\leq C \bigg(|\gamma| + \frac{(\ln n)^{1/p}}{n^{2(h(y_{\epsilon})-h_{\min})}}\bigg).
\end{equation*}
From the previous inequality, $ \lsp U_{n}(\gamma, \epsilon) \leq C |\gamma|$ and thus we get $\underset{\gamma\to 0}{\overline{\lim}}\  \lsp U_{n}(\gamma, \epsilon)=0$, which completes the proof.
\qed
\begin{rem} \rm
\label{ofjerifrofjerfioerj}
In the previous proof (in the case $q=0$), using
\eqref{zoiejozedijezoidjzeodijzeodijzzz1}, one can also see that the bias related to the convergence of  $\mu^{(2)}_{n}$ to $0$ is of order $1/\ln n$. \qed
\end{rem}

\subsection{Proof of Theorem \ref{th2} } 

In the first step we will find an upper bound for the covariance function of the increments of $B^h$. We
define
\[
r_n(i,j):= \text{cov} \left( B^h_{\frac{i+k}{n}} -  B^h_{\frac{i}{n}}, B^h_{\frac{j+k}{n}} -  B^h_{\frac{j}{n}} \right),
\qquad k=1,2.
\]
Recalling the notation at \eqref{cov}, we conclude the identity 
\[
r_n(i,j) = R_h\left(\frac{i+k}{n}, \frac{j+k}{n}\right) - R_h\left(\frac{i}{n}, \frac{j+k}{n}\right) 
- R_h\left(\frac{i+k}{n}, \frac{j}{n}\right) + R_h\left(\frac{i}{n}, \frac{j}{n}\right).
\] 
Since $h \in C^2([0,1])$ and the function $c$ defined at \eqref{cx} is a $C^{\infty}((0,1))$-function, 
we deduce by an application of Taylor expansion
\begin{align} \label{covest}
|r_n(i,j)| \leq n^{-2} \sum_{l,l'=1}^2 |\partial_{ll'} R_h(\psi_{ij}^n)| \qquad \text{for } |i-j|>2, 
\end{align}
where $\partial_{ll'} R_h$ denotes the second order derivative in the direction of $x_l$ and $x_{l'}$,
and $\psi_{ij}^n \in (i/n,(i+k)/n) \times (j/n,(j+k)/n)$. Now, we will compute an upper bound for the right side of \eqref{covest} for $i \not = j$. First, we observe that 
\[
R_h(t,s) = F(t,s)\ G(t,s, h(t)+h(s)),
\] 
where 
\begin{align*}
F(t,s) = \frac{ c^2_{h_{t,s}}}{c_{h(t)}c_{h(s)}}, \qquad
G(t,s, H) = \frac{1}{2} \left( {|t|}^{H} + {|s|}^{H}  - {|t-s|}^{H} \right).
\end{align*}
We remark that $G(t,s, 2H)$ is the covariance kernel of the fractional Brownian motion with Hurst parameter $H \in (0,1)$. 

Since $h \in C^2([0,1])$, $c \in C^{\infty}((0,1))$ and $c_x \not =0$ for $x \in (0,1)$, we conclude that 
\[
|\partial_l F(t,s)|, ~  |\partial_{ll'} F(t,s)| \leq C, \qquad l,l'=1,2, \quad (t,s)\in [0,1]^2.
\]
We concentrate on the second order derivative $\partial_{11} R_h(\psi_{ij}^n)$; the estimates for the other second order derivatives are obtained similarly. We have that 
\begin{align*}
\partial_{11} R_h(t,s) &= \partial_{11} F(t,s)\cdot G(t,s, h(t)+h(s)) \\[1.5 ex]
&+ 2\partial_{1} F(t,s)
\left[ \partial_{1} G(t,s, h(t)+h(s)) + h'(t) \cdot \partial_{3} G(t,s, h(t)+h(s))\right] \\[1.5 ex]
&+ F(t,s) \left[ \partial_{11} G(t,s, h(t)+h(s)) + 2  h'(t) \cdot \partial_{13} G(t,s, h(t)+h(s)) \right). \\[1.5 ex]
&\hspace{1.75cm}+\left.  h''(t) \cdot \partial_{3} G(t,s, h(t)+h(s)) + {(h'(t))}^{2}\cdot \partial_{33} G(t,s, h(t)+h(s)) )\right].
\end{align*}  
For the derivatives of the function $G$, we deduce the following estimates
\begin{align*}
|\partial_{1} G(t,s, h(t)+h(s))| &\leq C \left(t^{h(t)+h(s) -1} + |t-s|^{h(t)+h(s) -1} \right), \\[1.5 ex]
|\partial_{3} G(t,s, h(t)+h(s))| &\leq C \left(-\ln t\cdot t^{h(t)+h(s)} -  \ln s\cdot s^{h(t)+h(s)} 
-\ln |t-s|\cdot |t-s|^{h(t)+h(s)}  \right)  \\[1.5 ex]
|\partial_{11} G(t,s, h(t)+h(s))| &\leq C \left(t^{h(t)+h(s) -2} + |t-s|^{h(t)+h(s) -2} \right) \\[1.5 ex]
|\partial_{13} G(t,s, h(t)+h(s))| &\leq C \left( (1-\ln t) \ t^{h(t)+h(s) -1}
+(1-\ln |t-s|) |t-s|^{h(t)+h(s) -1} \right) \\[1.5 ex]
|\partial_{33} G(t,s, h(t)+h(s))| &\leq C  \left(\ln^2 t\cdot t^{h(t)+h(s)} +  \ln^2 s \cdot s^{h(t)+h(s)} 
+\ln^2 |t-s|\cdot |t-s|^{h(t)+h(s)}  \right),
\end{align*} 
which hold for $t,s \in (0,1]$ with $t \not = s$ and the third inequality holds whenever $h(t)+h(s) \not =1$ (if $h(t)+h(s)=0$ we simply have $\partial_{11} G(t,s, h(t)+h(s))=0$). Similar formulas and bounds are obtained for other second order derivatives of $R_h$. 
Using the boundedness of functions $F$, $h$ and its derivatives, together with the above estimates and  \eqref{covest} we obtain the inequality  
\begin{align} 
|r_n(i,j)| &\leq Cn^{-h(i/n)-h(j/n)} \left( i^{h(i/n)+h(j/n) -2} + j^{h(i/n)+h(j/n) -2}\right. \nonumber \\[1.5 ex]
\label{finalest} &\left. +|i-j|^{h(i/n)+h(j/n) -2} \right) \\[1.5 ex]
&\leq C n^{-2h_{\min}} \left( i^{2h_{\min} -2} + j^{2h_{\min} -2} +|i-j|^{2h_{\min} -2} \right)
, \qquad i,j \geq 1, |i-j| > 2. \nonumber
\end{align} 
When $|i-j|\leq 2$ we deduce from \cite[Lemma 1 p.$13$]{BS10} that
\begin{align} \label{finalest2}
|r_n(i,j)| \leq \text{var} \left( B^h_{\frac{i+k}{n}} -  B^h_{\frac{i}{n}}\right)
+ \text{var} \left( B^h_{\frac{j+k}{n}} -  B^h_{\frac{j}{n}}\right)\leq C n^{-2h_{\min}}.
\end{align}
We recall the identity $\text{cov}(Z_1^2,Z_2^2)= 2 \text{cov}(Z_1,Z_2)^2$ for a Gaussian vector $(Z_1, Z_2)$. By \eqref{finalest} and \eqref{finalest2}  we immediately conclude that 
\begin{align} \label{var}
\text{var}(V(B^h;k)^n) \leq C n^{-4h_{\min}+1} \sum_{i=1}^n i^{4h_{\min} -4} \leq C
\begin{cases}
n^{-4h_{\min}+1}  & h_{\min} \in (0,3/4) \\
\ln n\cdot n^{-2} & h_{\min} = 3/4 \\
n^{-2}  & h_{\min} \in (3/4,1)
\end{cases} 
\end{align}
Observing the decomposition 
\begin{align*}
S_n - \left( \frac12 \right)^{2h_{\min}} &= \frac{V(B^h;1)^n - \E[V(B^h;1)^n]}{V(B^h;2)^n} 
- \uh \frac{V(B^h;2)^n - \E[V(B^h;2)^n]}{V(B^h;2)^n} \\[1.5 ex]
&+ \left(\uh - \left( \frac 12 \right)^{2h_{\min}}\right)
\end{align*}
and in view of Proposition \ref{firstresult}, it is sufficient to show that 
\begin{align} \label{ratioconv}
\frac{\sqrt{\text{var}(V(B^h;l)^n)}}{\E[V(B^h;k)^n]} \to 0, \qquad k,l=1,2
\end{align}
to prove Theorem \ref{th2}. We assume again without loss of generality that $q=0$, $m=1$ and 
${\cM}_{h}=h^{-1}\{h_{\min}\}=\{x\}$. Using the notations  from the previous subsection together with the inequalities \eqref{oifjefoijefijoreijr1dede}
and \eqref{oifjefoijefijoreijr2dede}, we deduce the following lower bound, for $n$ large enough and for $\epsilon$ small enough:
\begin{align*}
\E[V(B^h;k)^n] & \geq \underline{\Gamma}^{(1)}_{n,k}(\gamma, \epsilon) 
+ \underline{\Gamma}^{(2)}_{n,k}(\gamma, \epsilon) \geq C_{\epsilon} \frac{n^{1-2h_{\min}}}{ 
(\ln n)^{1/p}}.  
\end{align*}
Thus, in view of  \eqref{var} we readily deduce the convergence at 
\eqref{ratioconv} for any $h_{\min} \in (0,1)$, which completes the proof of Theorem \ref{th2}. \qed


\begin{small}
\bibliographystyle{alpha}

\end{small}

\end{document}